\input amstex
\documentstyle{amsppt}
\input epsf.tex
%
%

\def\Int{\mathop{\roman{Int}}\nolimits}

\def\conj{\mathop{\roman{conj}}\nolimits}

\def\d{\partial}

%
%

\def\G{\Gamma}
\def\D{\Delta}
\def\e{\varepsilon}

\def\Con{\mathop\roman{Con}\nolimits}
%
%
\def\C{{\Bbb C}}
\def\R{{\Bbb R}}

%
%
\def\nologo{\let\logo@\relax}
%
%

\def\lle{<\!<}
\let\ge\geqslant
\let\le\leqslant

\def\barCP#1{\overline{\C\roman P}^{#1}}
\def\Cp#1{\C\roman P^{#1}}

\def\Rp#1{\R\roman P^{#1}}
%
%
\let\tm\proclaim
\let\endtm\endproclaim

\let\rk=\remark
\let\endrk=\endremark

%
%
%
%
\def\oo{\varnothing}
\def\A{\frak A}
\def\res{\roman{res}}

\magnification\magstep1
\nologo
\NoBlackBoxes
\NoRunningHeads

\topmatter
\title
A Rokhlin Conjecture and smooth quotients by the complex conjugation of
singular real algebraic surfaces
\endtitle
\author Sergey Finashin
\endauthor
\address
Middle East Technical University,
Ankara 06531 Turkey
\endaddress
\email serge\,\@\,rorqual.cc.metu.edu.tr
\endemail
\comment
\abstract
Suppose $\C X\to\Cp2$ is the double covering branched along a complex
curve $\C A\subset\Cp2$ defined over $\R$ and $\overline X=\C X/\conj$
is the quotient by a complex conjugation, $\conj$, lifted from $\Cp2$.
I prove that $\overline X\#_k\barCP2\conj\#_n\Cp2\#_m\barCP2$
for some $k,n,m\ge0$, if $\R A=\C A\cap\Rp2\ne\oo$ and
$\C A$ splits into a union $\C A=\C B\cup\C C$ of non-singular
transverse curves defined over $\R$.

I characterize the singularities of a complex surface, $\C X$,
defined over $\R$, whose quotient $\overline X=\C X/\conj$ is smooth.
\endabstract
\endcomment
\endtopmatter

\document
\heading
\S1 Introduction
\endheading

\subheading{1.1. A Rokhlin Conjecture and its variations}
Considering real algebraic varieties, I use prefixes $\C$ and $\R$
to denote their complex point sets and the real point sets
respectively, and put a bar to denote the orbit space for the
involution of the complex conjugation, $\conj$, for example, $\C
X$, $\R X$, $\overline X=\C X/\conj$. I identify  $\R X$ with its
image in $\overline X$. The same conventions will be used in the
notation of $\conj$-invariant subsets of a real algebraic variety
$\C X$.

If $\C X$ is a non-singular complex surface, then $\overline X$ is
a $4$-manifold and the quotient map $q\:\C X\to\overline X$ is a
double covering branched along $\R X$. One can endow $\overline X$
with an orientation and a smooth structure making $q$ smooth and
orientation preserving, for example $\overline{\roman
P}^2=\Cp2/\conj$ is well known to be diffeomorphic to $S^4$.

In 1980's I learned from V.A.Rokhlin several problems about
the topology of $\overline X$ for double planes $\C X$,
motivated by studying the topology of their branching loci,
the curves $\C A\subset\Cp2$.
One of these problems is as follows.
A curve $\C A\subset\Cp2$,
defined by a real form $f(x,y,z)$ of degree $d=2k$, with $\R A\ne\oo$,
divides $\Rp2$ into two regions,
$\Rp2_\pm=\{[x:y:z]\in\Rp2\,|\,\pm f\ge0\}$.
If $\R A$ is non-singular, then
one of these regions is orientable, the other is not, and
for the sake of definiteness it is usually assumed
that $\Rp2_+$ is the orientable region.

The unions $\A_\pm=\Rp2_\pm\cup\overline A$ are closed surfaces in
$S^4=\overline{\roman P}^2$, called by Rokhlin {\it the Arnold surfaces}.
They consist of the two smooth pieces, $\Rp2_\pm$ and $\overline
A$, which meet normally along $\R A$, and therefore $\A_\pm$ can
be smoothed along $\R A$. Rokhlin suggested to prove that the
embeddings $\A_\pm\subset S^4$ are {\it standard}, in the sense
that they can be obtained via ambient connected sum from a few
copies of a standard embedding of a torus, $T^2\subset S^4$, if
$\A_\pm$ is orientable, or from standard embeddings of $\Rp2$ if
$\A_\pm$ is not orientable (recall that there exist two isotopy
classes of standard embeddings of $\Rp2$ into $S^4$, which are
distinguished by the normal Euler number, $2$ or $-2$).

The surfaces $\A_\pm$ appeared in \cite{A1} in connection with
studying the topology of the double planes $p\:X\to\Cp2$ branched along
$\C A$. There exist $2$ liftings of the complex conjugation
from $\Cp2$ to $X$, and I denote them by $\conj^\pm$.
The double plane $X$ endowed with an involution
$\conj^\pm$ is denoted by $\C X^\pm$. The reason behind
 such a notation is that one can view $\C X^\pm$ as
a real algebraic surface defined
by the equation $w^2=\pm f(x,y,z)$ in a weighted projective $3$-space.
As follows from \cite{Ar}, the projection
$q^\pm\:\overline X^\pm\to S^4=\overline{\roman P}^2$
induced by $p$ is a double covering branched along $\A_{\mp}$.
Thus, if the Arnold surface $\A_\mp$ is standard, then $\overline X^\pm$
is decomposable into a connected sum of a few copies of $\Cp2$, $\barCP2$,
if $\A_\mp$ is non-orientable, or copies of $S^2\times S^2$ if orientable.
Such a splitting will be called a complete decomposition.
A weaker form of Rokhlin's question: {\it is it true that
$\overline X^\pm$ admits such a decomposition
for any double plane $\C X^\pm$ with $\R X^\pm\ne\oo$}.
I will refer to this statement as to
{\it CDQ-Conjecture}  (an abbreviation for
``Complete Decomposability of the Quotients'').
A review of the known result concerning the CDQ-Conjecture
and its natural extensions to the other real algebraic surfaces
can be found in \cite{F2}.

A weaker version of the CDQ-Conjecture was suggested by Akbulut \cite{Ak}:
{\it to prove that SW (Seiberg-Witten) invariants of $\overline X^\pm$
vanish if $b_2^+(\overline X^\pm)>1$}.
As an intermediate variation, I present below one more conjecture.
Let us call 4-manifolds $X$ and $Y$ {\it blow-up stable equivalent}
(BUS-equivalent) if
$X\#_n\barCP2$ is diffeomorphic to $Y\#_m\barCP2$ for some $n,m\ge0$.
If $X$ is BUS-equivalent to $\#_n\Cp2$, $n\ge0$, then we call
$X$ {\it BUS-trivial}.
Complete decomposability implies BUS-triviality, but not vice versa
(some algebraic surfaces with the inverted orientation
are BUS-trivial, as follows from \cite{MM}, but not decomposable).
On the other hand, BUS-triviality guarantees vanishing of the Donaldson
and SW-invariants
(more generally, vanishing of these invariants is a property which is
preserved under BUS-equivalence, cf. \cite{FS}).
So, {\it the BUS-Conjecture} suggesting that
the quotients $\overline X^\pm$ are BUS-trivial for the double planes,
$\C X^\pm$ with $\R X^\pm\ne\oo$, is weaker then CDQ-conjecture, but
stronger then Akbulut's conjecture.

\subheading{1.2. The results}
The aim of this paper is to discuss an extension of the Rokhlin Conjecture
 and its versions to the case of real surfaces with singularities
and to present some new results in its support.
The main goal is the following theorem.

\tm{1.2.1. Theorem}
Assume that a real plane curve, $\C A_0$, of degree $2k$ and with
$\R A_0\ne\oo$, splits into a union
$\C A_0=\C B_0\cup \C C_0$ of transverse non-singular curves and
a non-singular curve, $\C A$, is obtained from  $\C A_0$ by a small
perturbation. Let $\C X^\pm$ denote the double plane branched along $\C A$.
Then the both quotients,  $\overline X^+$ and  $\overline X^-$, are
BUS-trivial.
In particular, the SW-invariants of $\overline X^\pm$ vanish, if $k>3$.
\endtm

The condition $k>3$ is required for
$b_2^+(\overline X^\pm_t)=\frac12(k-1)(k-2)>1$.
The proof of Theorem 1.2.1, which is given in \S3,
is based on the analysis of the bifurcations
which can experience $\overline X^\pm$ as we deform $\C B_0$ and $\C C_0$.
One more condition of vanishing of SW-invariants
is given by Theorem 2.3.1.

\heading \S2. Rokhlin's Conjecture and its versions in the case of
singular real algebraic surfaces
\endheading

\subheading{2.1. Smoothly-folding real surface singularities}
One or the both of the Arnold surfaces, $\A_\pm$, may turn out to be
smoothable even if we admit certain types of real singularities on
$\C A$. The condition of smoothability is that
the links of $\A_\pm\subset S^4$ at the singular points,
$x\in\R A$, are unknots. These unknots are piecewise smooth and their
smoothings can be easily extended to
smoothings of the embeddings $\A_\pm\subset S^4$.
This gives a smoothing of the double covering
$\overline X^\pm\to S^4$, branched along $\A_\pm$, for
the double planes, $\C X^\pm\to\Cp2$, branched along $\C A$.

A singularity at a point $x\in\R X$ of a real algebraic surface
will be called
{\it smoothly-folding}, or briefly, SF-singularity, if the link of
$\overline X$ at $x$ is a $3$-sphere (so, $\overline X$ is a manifold
around $x$).
We call a real algebraic surface $\C X$ {\it SF-surface} if it has only
SF-singularities.
Like in the case of the double planes, a smoothing of $\overline X$
at the singularities can be easily extended to a smoothing
along the whole $\R X$.

{\it A cone-like node} (like the locus $\{x^2+y^2-z^2=0\}\subset\C^3$)
is an example of SF-singularity, whereas {\it a dot-like node}
(like the locus $\{x^2+y^2+z^2=0\}$) is not.
Another example of an SF-singularity appears on the double plane
$\C X^\pm$, if the branching locus, $\C A$, splits at $x\in\R A$
into $m$ non-singular real and pairwise transverse branches, where $m$
is odd.
More examples: any simple surface singularity has a real form
which is an SF-singularity (I mark this form by a superscript
``$-$'', e.g., $A_n^-$, $D_n^-$, $E_n^-$).

More generally, SF-singularities can be characterized
in terms of a minimal resolution,  $\res\:\C X^\res\to\C X$,
and the exceptional curve, $\C E\subset\C X^\res$, over $x$
(a minimal resolution of $\C X$ yields a real algebraic surface,
because it is obtained after several blows up at a real point
or at a pair of $\conj$-symmetric imaginary points of a real surface;
for more details about real resolutions see, e.g., \cite{S}).
Since $\res$ may be not {\it a good} resolution, we have to modify
the usual definition of the resolution graph as follows.
The vertices of a graph $\G$ split into 2 types: the vertices
of the first type are in $1$--$1$ correspondence with the irreducible
components of $\C E$ and the vertices of the other type are in
$1$--$1$ correspondence with the intersection points of these components.
A pair of vertices of $\G$ is connected by an edge if one of the vertices
represents a component of $\C E$ and the other represents a point
on this component.
The complex conjugation gives an involution
on $\G$ and we denote by $\overline{\G}$
the quotient graph with respect to this involution.

\tm{2.1.1. Theorem}
Assume that $x\in\R X$ is an isolated normal real surface singularity,
and $\C E$ is the exceptional curve of its minimal (real) resolution.
Then the following conditions are equivalent.
\roster
\item
$x$ is an SF-singularity,
\item
$\overline E$ is contractible,
\item
$\C E$ splits into real rational irreducible components,
$\C E=\C E_1\cup\dots \C E_n$, with $\overline E_i$, $i=1,\dots,m$,
being homeomorphic to a $2$-disc $D^2$,
and with the quotient graph $\overline{\G}$ being a tree.
\endroster
\endtm

\rk{Remark}
If a minimal resolution turns out to be {\it good}, then we may use
the usual resolution graph instead of $\G$ (in this case,
the graph $\G$ is just a subdivision of the resolution graph).

Note also that the condition $\overline E_i\cong D^2$ means that $\C E_i$
is a topological sphere with $\R E_i\ne\oo$, or is obtained from
such a sphere after identifying a few pairs of $\conj$-symmetric points.
\endrk

\demo{Proof}
The equivalence (2)$\Leftrightarrow$(3) is trivial.
Let us assume (2) and
consider a compact $\conj$-invariant cone-like neighborhood,
$\C U\subset\C X$ of $x$. Let $\C U^\res\to\C U$ denote
the restriction of the resolution $\res\:\C X^\res\to\C X$ over $\C U$.
Note that $\C U^\res$ is a regular neighborhood of $\C E$ in $\C X^\res$
(i.e., $\C E$ is a spine of $\C U^\res$).
This implies that $\overline U^\res$ is a regular neighborhood of
$\overline E$, and thus, of the graph $\overline{\G}$
(which is a simple deformational retract of $\overline E$).
Therefore, $\overline U^\res\cong D^4$ and
$\d (\overline U)=\d (\overline U^\tau)\cong S^3$.

Assuming (1), we can easily obtain, using the exact homotopy
(or homology) sequence of the pair
$(\overline U^\res,\d(\overline U^\res))$, that $\overline U^\res$,
and therefore its spine $\overline E$, is simply connected;
here we do not even use that the resolution, $\res$, is
{\it minimal}.
It follows that $\overline \G$ is a tree and that
all the irreducible components of $\C E$ are rational
and their quotients in $\overline E$ can be homeomorphic
either to $D^2$ (if a component is real, i.e., $\conj$-invariant)
or to $S^2$ (if it is not real).
To complete the proof, we need to show
that all the components of $\C E$ are real.

Blowing up $\C U^\res$, we can obtain a real {\it good resolution},
$\C U^\res_g\to \C U$, with the exceptional
curve, $\C E_g\subset \C U^\res_g$, consisting of
non-singular irreducible components, which intersect pairwise
transversally and do not have triple intersection points.
We can associate graphs
$\G_g$, $\overline{\G}_g$, to $\C E_g$ like $\G$ and
 $\overline{\G}$ were associated to $\C E$.

Let us remove from $\overline{\G}_g$ the vertices corresponding to
the real (i.e., $\conj$-invariant) components of $\C E_g$
and to the real intersection points on these components,
together with the adjacent edges.
The rest splits into several
components, $C_1,\dots,C_n$, each component being a tree. The vertices of
the first type of
$C_i$ represent a set of transversal spheres embedded into
$\overline U^\res_g$.
A regular neighborhood of their union is a plumbing manifold, which will
be denoted by $V_i$.
Note that the intersection form in $V_i$
is negative, since the form in $\C U^\res_g$
(and, thus, in $\overline U^\res_g$) is negative.
 On the other hand, it is not difficult to see that $\overline U^\res_g$
is homeomorphic to a boundary-connected-sum of $V_1,\dots,V_n$,
which implies that $\d V_i\cong S^3$, $i=1,\dots,n$, since $\d
(\overline U^\res)\cong S^3$. Now we apply Neumann's
result \cite{N} stating that the boundary of a plumbing
manifold with a negative intersection form determines this
manifold up to blowing up and blowing down. This implies that all
the irreducible imaginary components of $\C E_g$ can be blown down
and do not appear in $\C E$.
\qed\enddemo

It is trivial that a blow-up at a non-singular point of $\R X$
does not change $\overline X$ (the proof:
$\overline X'\cong(\C X\#\barCP2)/\conj
\cong\overline X\#S^4\cong\overline X$).
The following observation generalizes this fact.

\tm{2.1.2. Corollary} Assume that $\C X'\to \C X$ is obtained by
blowing up at an SF-singularity. Then $\overline X'\cong\overline
X$.
\endtm

The proof is not a difficult exercise, after Theorem 2.1.1.
\qed

\subheading{2.2. The case of curves $\C A$ splitting into real lines}
In \cite{F1} I proved the Rokhlin Conjecture in the following special
cases.

\tm{2.2.1. Theorem. \cite{F1, Theorem 3.1}}
Assume that a curve $\C A_0\subset\Cp2$ splits into a generic arrangement
of $2k$ real lines, $k\ge1$ and $\C A$ is obtained by a real perturbation
of $\C A_0$.
Then the both Arnold surfaces, $\A^\pm$, are standard.
In particular, $\overline X^\pm$ is completely decomposable for
the double planes, $\C X^\pm\to\Cp2$, branched along $\C A$.
\endtm

The Rokhlin Conjecture holds for
the Arnold surfaces, $\A_0^\pm$ of the curve $\C A_0$, as well,
and the above theorem indeed implies it, because
for a certain canonical perturbation scheme
the curve $\C A$ have the Arnold surface $\A^\pm$ isotopic to $\A_0^\pm$
(such a canonical perturbation scheme for $\A^+$ is the opposite
to a canonical perturbation scheme for $\A^-$).
A change in the canonical perturbation scheme at any node of $\C A_0$,
effects to $\A^\pm$ as a surgery, which replaces a small $2$-disc around
the nodal point by a M\"obius band, which effects to $\overline X^\mp$
as a blow up (here  $\C X^\mp$ is, as usual, the corresponding double plane).

The arguments in \cite{F1} are applicable also for real curves on
a quadric, $\C A\subset\Cp1\times\Cp1$, of bi-degree $(2k,2l)$,
which are obtained by a perturbation from a singular curve, $\C
A_0$, splitting into $2k+2l$ generating lines of the quadric. By
blowing up at a real point and the further blowings down, the
quadric $\Cp1\times\Cp1$ is transformed into $\Cp2$ and the curve
$\C A_0$ becomes a configuration of $2k+2l$ real lines, $2k$ of
which pass through one point and $2l$ pass through another one.

The method of \cite{F1} is applicable, more generally, for curves
$\C A$ splitting into an arbitrary arrangement of $2k$ distinct
real lines. Let $\C\roman{Q}\to\Cp2$ be obtained by blowing up at
the singular points of $\C A$ whose multiplicity is even and
greater then $2$, and $\C B\subset\C\roman{Q}$ the proper image of
$\C A$. Note that $\overline{\roman{Q}}\cong S^4$, since a blow up
at a real point does not change the quotient. $\R\roman{B}$
divides $\R\roman{Q}$ into a pair of regions, $\R\roman{Q}_\pm$,
like in the case of $\Rp2$, and we can similarly define a pair of
the Arnold surfaces $\A_\pm=\overline B\cup\R\roman{Q}_\pm\subset
S^4$. Furthermore, we can similarly define the double covering,
$\C X^\pm\to\C Q$, branched along $\C B$, endowed with one of the
two liftings of the complex conjugations on $\C Q$; the induced
projection $\overline X^\pm\to\overline Q=S^4$ is again a double
covering branched along $\A_\mp$. Note that $\A_\pm$ is connected,
unless $\C A$ is a pencil, i.e., unless all the $2k$ lines have a
common point. In the case of a pencil, $\A_\pm$ is a union of
$2$-spheres which bound disjoint $3$-balls in
$\overline{\roman{Q}}$ (this follows from \cite{F1,Theorem 5.2}),
and, thus, $\overline X^\pm$ is diffeomorphic to a connected sum
$\#_{k-1}(S^1\times S^3)$.

\tm{2.2.2 Theorem}
Assume that $\C A\subset\Cp2$ splits into a configuration of $2k$ real
lines different from a pencil. Then $\A_+$ and  $\A_-$
are standard surfaces in $S^4=\overline{\roman{Q}}$ and therefore
$\overline X^+$ and $\overline X^-$ are completely decomposable.
\endtm

The proof is analogous to the proof of Theorem 3.1 in \cite{F1}.

\subheading{2.3. Vanishing of SW-invariants for $\overline X$ after
a perturbation of an SF-singularity}
\vskip0mm
\noindent
Assume that $\C X$ is a real surface having one singular point,
$x\in\R X$, and a non-singular surface $\C X'$ is obtained from $\C X$
by a real perturbation.
More precisely, we assume that
there exists an equivariant (with respect to the complex conjugation)
diffeomorphism between the complement $\C X-\Int(\C U)$ for
a compact regular cone-like neighborhood, $\C U\subset\C X$, of $x$,
and the complement $\C X'-\Int(\C U')$ of some compact $\conj$-invariant
codimension $0$ submanifold $\C U'\subset\C X'$ of a smooth
real algebraic surface $\C X'$.
A standard example that I mean is $\C U'$ being a Milnor fiber
in the case of a complete intersection singularity
(for instance, $\C X$ may be a double plane branched along a curve,
$\C A$, with one singular point, and $\C X'$ a double plane branched
along a non-singular curve $\C A'$ obtained by a perturbation of $\C A$).

\tm{2.3.1. Theorem}
Assume that $\C X$ and $\C X'$ are as above, $x\in\R X$ is
an  SF-singularity and $0<p_g(\C X^\res)<p_g(\C X')$, where
$p_g$ denotes the geometric genus and $\C X^\res\to\C X$ a real resolution.
Then SW-invariants of $\overline X'$ vanish.
\endtm

\demo{Proof} $\d (\overline U)\cong S^3$, for SF-singularities,
thus, $\overline X'$ is obtained by taking a connected sum of
$\overline X$ with a 4-manifold $\hat U'$, obtained by attaching a
$4$-ball to $\overline U'$.

Standard and well known calculations imply that
$b_2^+(\overline X')=p_g(\C X')$ and $b_2^+(\overline X)=p_g(\C X^\res)$
(see for instance Lemma 4.2.2 in \cite{F3}),
so, $\overline X'$ have trivial SW invariants,
as a connected sum of manifolds with non-negative intersection forms.
\qed
\enddemo

\rk{Remark}
The above arguments can be obviously applied as well to
SF-surfaces with several singular points, provided the same inequalities
hold.

Note also that there exist vanishing theorems for SW-invariants,
which generalize the connected sum theorem that was used. Namely, such
theorems concern $4$-manifolds which can be split
along some other types of $3$-manifolds (rather then along
$S^3$, like in the case of the connected sums).
Accordingly, Theorem 2.3.1 can be modified and extended to the corresponding
class of singularities.
\endrk

\subheading{2.4. The singular and the local versions of the Rokhlin problem}
A connected sum splitting of $\overline X'$ in
the above proof of Theorem 2.3.1 gives rise to the questions
about the topology of $\C X$ and $\C U'$; for instance, is it true that
$\overline X$ and $\hat U'$ are always completely decomposable whenever
they are simply connected ?

\rk{Example} Assume that $\C A$ is {\it almost a pencil}, that is,
a union of $2k$ real lines, which all, except one, contain a
common point. Then $\A_\pm$ is an unknotted sphere and $\overline
X^\pm$ is diffeomorphic to $S^4$. This shows equivalence of the
``local version'' of the CDQ-conjecture with the usual ``global
version'', in the particular case of the real double planes, whose
branching locus is obtained from $\C A$ by a perturbation.
\endrk

\heading
\S3. Proof of Theorem 1.2.1
\endheading

\subheading{3.1. Simple nodal deformations of the real plane curves and the
double planes}
\newline
Recall that the real plane curves of degree $d$ constitute a real projective
space, $\Cal C_d$, of dimension $\binom{d+2}2$, where the singular curves
form a hypersurface, $\D\subset\Cal C_d$, called the discriminant
hypersurface.
  By a real deformation of a curve I mean a continuous family of curves,
 $\C A_t\subset\Cp2$, $t\in[a,b]$, that is a path in $\Cal C_d$.
If a path is generic, it intersects $\D$ transversally at non-singular
points, which are represented by the curves having a nodal singularity
(an ordinary double point).
In this case, $\C A_t$,  $t\in[a,b]$, is called {\it a simple
nodal deformation}.

There may be two types of nodes on real curves:
dot-like nodes, like the locus $\{x^2+y^2=0\}\subset\C^2$, and cross-like
nodes, like $\{x^2-y^2=0\}$.
Assume that $d$ is even and $\C A_t$, $t\in[0,1]$, is a simple  nodal
deformation. Taking the double planes branched along $\C A_t$, we obtain
two families, $\C X_t^\pm$.
Let us fix one of these two continuous families and denote it by
$\C X_t$,
omitting the superscript in the notation to avoid an ambiguity, which
appears if $\C X_0^+$ is continuously transformed into $\C X_1^-$
(this may happen because the
sign-superscripts are not defined canonically for nodal curves,
like for non-singular ones).
Furthermore, I denote by $W_t$ the corresponding region, $\Rp2_+$
or $\Rp2_-$, namely, $W_t=p_t(\R X_t)$, where $p_t\:\C X_t\to\Cp2$
is the branched covering.
The choice of $\C X_t$ is obviously determined by the choice of $W_0$.

The cross-like nodes of $\R A_t$ are covered obviously by the cone-like
nodes of $\R X_t$.
A dot-like node at $x\in\R A_t$ is covered by a cone-like node of $\R X_t$
provided $x$ lies in the interior of $W_t$, otherwise it is covered
by a dot-like node of $\R X_t$.
If in a process of deformation a cone-like node appears on
$\R X_{t}$, then one of the quotients,  $\overline X_{t+\e}$ or
 $\overline X_{t-\e}$, for $0<\e\lle1$, is diffeomorphic to $\overline X_t$
and the other one differs by a blow-up, i.e., is diffeomorphic to
$\overline X_t\#\barCP2$.
In the case $\R X_t$ has a dot-like node, $\overline X_t$ is not a manifold
and $\overline X_{t+\e}$ differs from $\overline X_{t-\e}$ by
a rational blow-down (or a rational blow-up) of multiplicity $2$,
in the sense of \cite{FS} (see \cite{Le}, \cite{Ak} or \cite{F1}).
So, the difficulty to prove that for {\it an arbitrary} non-singular
double plane, $\C X$, the quotient $\overline X$ is BUS-trivial consists in
proving that its branching locus, $\C A\subset\Cp2$, can be connected
by a {\it nice} nodal deformation
with some other curve, $\C A'$, with a BUS-trivial quotient $\overline X'$
of the corresponding double plane $\C X'$. Here, a nodal deformation is
{\it nice} if dot-like nodes do not appear on $\R X_t$.

\subheading{3.2. Simultaneous nodal deformations} One can handle
with the above difficulty if a curve $\C A$ is obtained after a
non-singular perturbation of a curve,
 $\C A_0$, splitting into a union $\C B_0\cup \C C_0$ of two transverse
non-singular real curves. In this case, we consider a deformation
$\C A_t$, $t\in[0,1]$, which splits as $\C A_t=\C B_t\cup\C C_t$,
where $\C B_t$ and $\C C_t$ are simple nodal deformations. In
addition to the nodes, we allow $\C A_t$ to have real
singularities of the type $A_3$, which appear at the points of a
simple (quadratic) tangency of $\R B_t$ and $\R C_t$. In a generic
pair of nodal deformations, no other singularities are possible
(so, none of the curves, $\C B_t$, $\C C_t$, passes through the
singularities of the other one, there is no higher order tangency
points and no imaginary tangency points). Such a deformation,  $\C
A_t$, is to be called {\it a simultaneous nodal deformation}. If
moreover, the fixed covering family of the double planes, $\C
X_t$, contains no dot-like nodes, then we call  $\C A_t$, {\it a
nice simultaneous nodal deformation}.

Given such a deformation, the quotients $\overline X_t$, $t\in[0,1]$, appear
to be smooth $4$-manifolds in the complement of a few singular points.
A singular point, $x\in\overline X_t$, may appear either from a conjugated
pair of imaginary intersection points of $\C B_t$ and $\C C_t$, or from
a tangency point of $\R B_t$ and $\R C_t$.
In the first case, $x$ is obviously again a nodal singularity.
In the second case, a singularity of $\C X_t$ of the complex type $A_3$
may have two real forms: $A_3^+$, or $A_3^-$,
the first of which can be represented by the model
$\{x^4-y^2=z^2\}\subset\C^3$, and the second by the model
$\{x^4-y^2=-z^2\}$.
If $x\in\C X_t$ has type $A_3^-$, then $\overline X_t$ is smoothed
at $x$ (since $A_n^-$ is an SF-singularity).
If $x\in\C X_t$ has type $A_3^+$, then
the link of $\overline X_t$ at $x$ is $\Rp3$,
so $\overline X_t$ has at $x$ a nodal singularity.

\subheading{3.3. Smooth nodal $4$-manifolds}
By a smooth manifold with isolated singularities I mean
a topological space, $Y$, which is a smooth manifold outside
the singular locus, $\Sigma\subset Y$, constituted by
a discrete set of points, each of which has a smooth cone-like neighborhood.
Such a neighborhood around  $y\in\Sigma$ is by definition a compact
neighborhood, $U_y\subset Y$, endowed with a homeomorphism,
$\phi_y\:U_y\to\Con(L_y)$, which is a diffeomorphism outside $y$.
Here $L_y$ is a $3$-manifold called {\it the link of $Y$ at $y$}
and $\Con(L_y)$ is a cone with the base $L_y$ (i.e., a join of
$L_y$ and $y$).
We may view $(U_y,\phi_y)$ as a ``chart'' required to bring a
``smooth'' structure to $y$.
It is well known that algebraic varieties with isolated singularities
have such a structure. The quotients $\overline X$ of real algebraic
surfaces,
$\C X$, having only isolated singularities, give another example.
We consider below only a special case of smooth $4$-manifolds
with all the links $L_y$ being homeomorphic to $\Rp3$; such manifolds
will be called {\it smooth nodal $4$-manifold}.

Given such a manifold, $Y$, we define a {\it resolution}
of a singularity at $y\in Y$, as a surgery replacing a cone-like
neighborhood, $U_y$,
by the total space of a smooth $D^2$-fiber bundle over $S^2$ with the
normal Euler number $-2$ (that is the usual surgery characterizing
a resolution of an algebraic node).
Let us denote by $Y^\res$ the manifold obtained after resolution of
all the nodes of $Y$. The differential type of $Y^\res$ is
well defined, because the orientation-preserving diffeomorphism
group of $\Rp3$ is known to be connected (cf. \cite{H}).
We call a pair of smooth nodal 4-manifolds,
$Y_1$, $Y_2$, BUS-equivalent, if $Y^\res_1$ is BUS-equivalent to
$Y_2^\res$. We call such a nodal manifold, $Y$, BUS-trivial,
if $Y^\res$ is BUS-trivial.

\subheading{3.4. Proof of Theorem 1.2.1} Theorem 1.2.1 follows
from the lemmas below.

\tm{3.4.1. Lemma}
Assume that $\C A_t=\C B_t\cup\C C_t$, $t\in[0,1]$, is a nice simultaneous
nodal deformation of plane real curves, of degree $\deg(\C A_t)=2k$, $k\ge1$,
and $\C X_t$ the associated deformation of the double planes.
Then $\overline X_0$ is BUS-equivalent to $\overline X_1$.
\endtm

Here, as before, $\C X_t$ is the corresponding family of the double
planes, which is determined by fixing a choice of $W_0$.

\tm{3.4.2. Lemma}
Assume that plane real curves $\C A_0$ and $\C A_1$ of degree $2k$, $k\ge1$,
with $\R A_0\ne\oo$,  $\R A_1\ne\oo$,
split into unions
$\C A_0=\C B_0\cup\C C_0$ and $\C A_1=\C B_1\cup\C C_1$
of nonsingular transverse curves, so that $\C B_0$ and $\C B_1$ have the same
degree. Then $\C A_0$ and $\C A_1$ can be connected by
a nice (with respect to any fixed choice of $W_0$)
simultaneous nodal deformation.
\endtm

In combination with Theorem 2.2.1 this gives

\tm{3.4.3. Corollary}
Assume that a curve $\C A_0$ splits like in Lemma 3.4.2. above.
Then $\overline X_0$ is BUS-trivial.
\endtm

Perturbing the imaginary singularities of $\C A_0$ we resolve
the nodal singularities of $\overline X_0$; perturbing the real nodes
we blow up $\overline X_0$ or preserve it topologically, as was mentioned
in 3.1. So, after perturbation of all the nodes of $\C A_0$ we obtain
a BUS-equivalent (thus, a BUS-trivial) manifold. This
completes the proof of Theorem 1.2.1.
\qed

\rk{Remark}
The condition that the curves $\C B_0$ and $\C C_0$ are transverse
(like a generic position condition that the tangency points of
$\C B_t$ and $\C C_t$ should not have order higher then $2$)
is not very essential in Theorem 1.2.1.
This follows from analysis of deformations of a singularity at a point
$s\in \C A_0$, at which $\C B_0$ and $\C C_0$ have an oder $n$ tangency
(that is a simple singularity of the type $A_{2n-1}$).
More generally, using the well-known classification of deformations of
the real simple singularities (cf. \cite{Ch}), one can prove the following

\tm{3.4.4. Theorem} If a real surface $\C X$ has a simple
singularity at $s\in\R X$ which is not equivalent to the
singularity defined in $\C^3$ by the equation $x^{2n}+y^2+z^2=0$,
and $\C X'$ is obtained by a real non-singular perturbation of $\C
X$, then the BUS-equivalence class of $\overline X'$ is
independent of the perturbation.

For the above exceptional singularity, the same claim is true for
all the perturbations except the one for which the real locus, $\R X'$,
vanishes near $s$.
\endtm
\endrk

In fact, $\overline X$ has at $s$ a singularity topologically
equivalent to an algebraic surface singularity and $\overline X'$
is BUS-equivalent to $\overline X^\res$, obtained from $\overline
X$ by resolution of this singularity. There exists a conceptual
proof of it, which explains this phenomenon better then a
straightforward routine analysis. It uses a version of the trick
applied by Donaldson \cite{D} to real K3 surfaces. This trick
consists in variation of the complex structure (within a family
parameterized by a twistor line), so that the involution of
complex conjugation becomes holomorphic. One can apply this trick
similarly for real simple singularities and their real
deformations, using the Kronheimer Torelli-like theorem \cite{Kr}
instead of Yau's theorem used by Donaldson.

\subheading{3.5. Proof of Lemma 3.4.1}
It is enough to determine the bifurcation of $\overline X_t$, after
passing a singularity of the type $A_3^\pm$
(see Figure 1).

\midinsert \vskip-5mm \botcaption{Figure 1} The singularities
$A_3^+$ and $A_3^-$. The region $W_t$ is shaded.
\endcaption
\vskip5mm \epsfxsize=13cm \centerline{\epsfbox{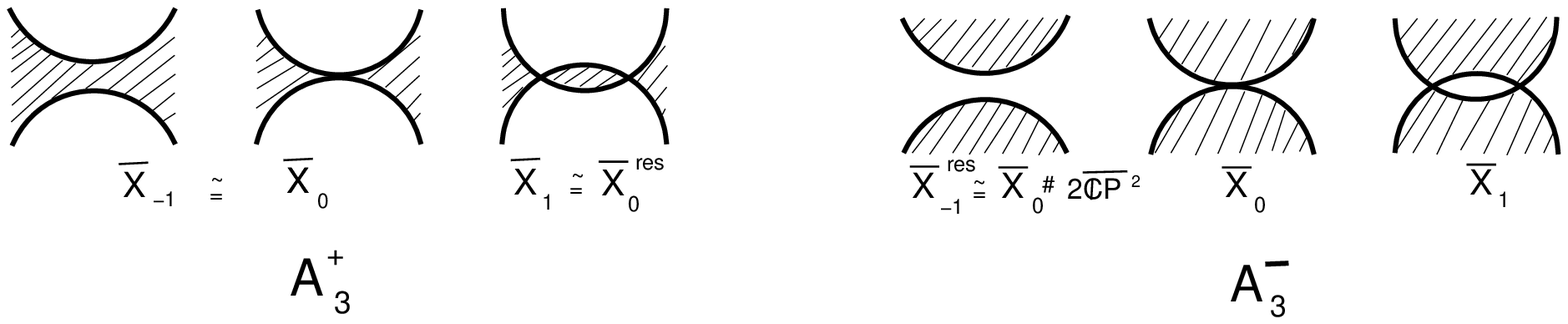}}
\endinsert

Since $\overline X_t\to S^4$ is a double covering branched along
the corresponding Arnold surface, $\A_t$, for $\C A_t$, the problem
is reduced to analysis of the bifurcations of $\A_t$
in the following model example.

Let $\C A_t\subset\C^2$, $t\in[-1,1]$, denote the zero locus
$\{f_t(x,y)=0\}$, where $f_t(x,y)=(x^2-t)^2-y^2$.
The family of surfaces $\C X_t\subset\C^3$ defined by the equation
$f_t(x,y)=z^2$ represents a deformation of a singularity $A_3^+$.
In this case, the quotient $\overline X_t\cong\overline X_0$ has one node
for $t\le0$, and $\overline X_t\cong\overline X_0^\res$, for $t>0$.
A deformation of a singularity $A_3^-$ is represented by a family
of surfaces $\C X_t=\{f_t(x,y)=-z^2\}\subset\C^3$.
In this case, $\overline X_t\cong\overline X_0$ is smooth, if $t\ge0$.
If $t<0$, then $\overline X_t$ has a node, such that
$\overline X_t^\res\cong\overline X_0\#2\barCP2$.
Thus, in the both cases,
$\overline X_{-1}$ is BUS-equivalent to $\overline X_1$.
\qed

\subheading{3.6. Proof of Lemma 3.4.2}
We can make ``nice'' a given simultaneous nodal deformation,
$\C A_t=\C B_t\cup\C C_t$ using the following trick.
Assume that $\C A_{t_0}$ has
a ``bad'' dot-like node, $x\notin\Int(W_{t_0})$, and
for definiteness, $x\in\R B_{t_0}$.
Let us assume $\C C_{t_0}$ is non-singular (just a generic position
assumption)
and that $\R C_{t_0}\ne\oo$ (a more essential assumption).
Choose a sufficiently small $\e>0$, so that
neither $\C B_{t}$, nor $\C C_{t}$ have singularities for
$|t-t_0|\le\e$, except the one at $x$ for $t=t_0$, and furthermore,
so that $x\notin\R C_t$ for $|t-t_0|\le\e$.
We may choose a loop, $[t_0-\e,t_0+\e]\to PGL(3;\R)$, $t\mapsto T_t$,
centered at the unity, and replace $\C C_t$, $t\in[t_0-\e,t_0+\e]$,
by a curve $T_t(\C C_t)$ obtained from $\C C_t$ by the action of $T_t$
(i.e., just move $\C C_t$ by linear projective transformations),
so that the real part $\R C_t$ passes  once across the point $x\in\Rp2$
as $t$ vary from $t_0-\e$ to $t_0$, see Figure 2.

\midinsert
\vskip-5mm
\botcaption{Figure 2}
Elimination of a ``bad'' dot-like
nodal singularity ($W_t$ is shaded)
\endcaption
\vskip5mm
\epsfxsize=13cm \centerline{ \epsfbox{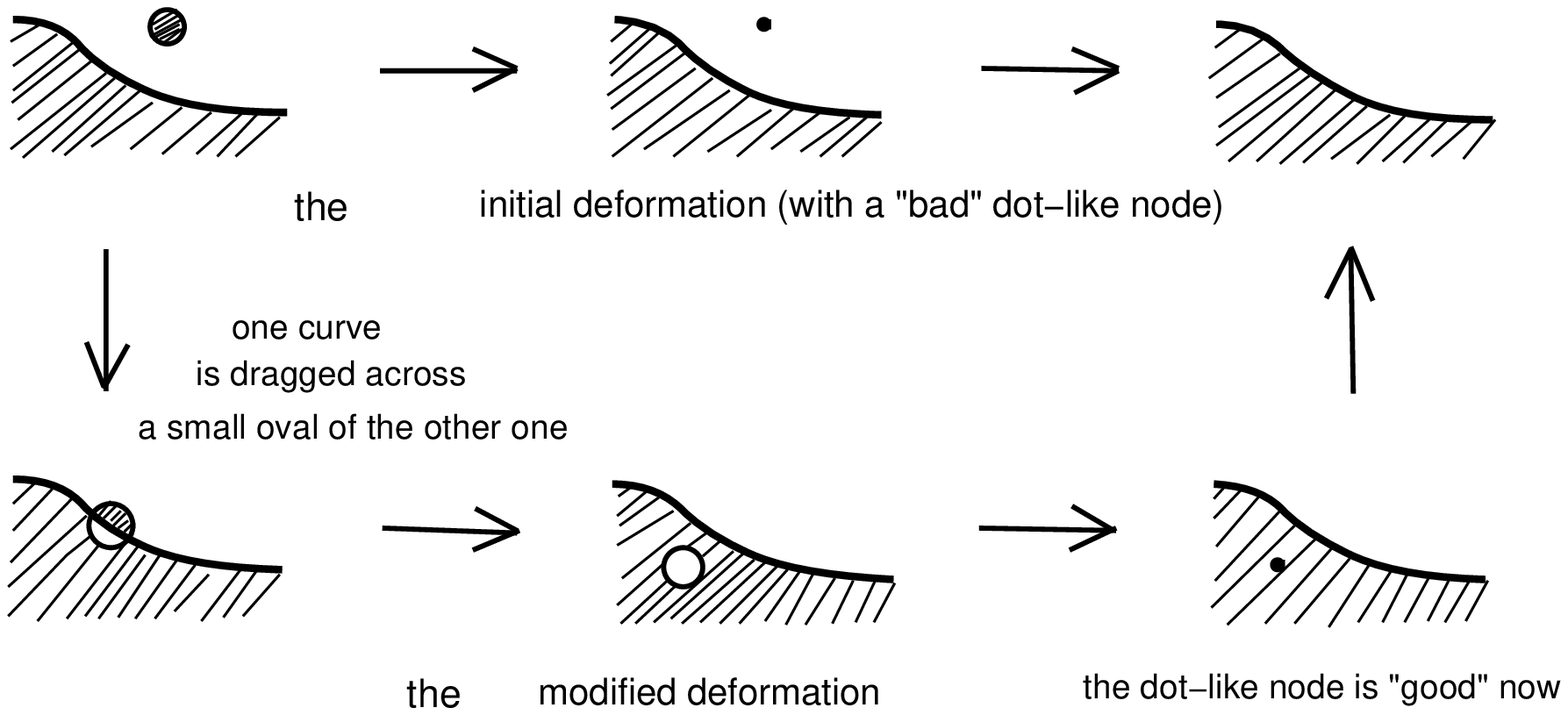}}
\endinsert
For a generic choice of the family $T_t$, this modification gives again
a simultaneous nodal deformation, but now
a dot-like node at $x$ appears in the interior of the (modified)
region $W_{t_0}$ (see Figure 2).

To apply the above trick and get rid of all the ``bad'' nodes, it
is required that $\C C_t\ne\oo$ if a ``bad'' node appears on $\C
B_t$ and $\C B_t\ne\oo$ if such a node appears on $\C C_t$. It is
not difficult to provide this as follows. Assuming that $\R
B_0\ne\oo$ and $\R C_1\ne\oo$, we deform first $\C C_0$ to obtain
$\C C_1$, not varying $\C B_0$ (allowing though linear projective
transformations, if needed to guarantee a generic position with
$\C C_1$), and then deform $\C B_0$ to obtain $\C B_1$. The case
$\R C_0\ne\oo$ and $\R B_1\ne\oo$ is analogous. If $\R C_0=\R
C_1=\oo$, but $\R B_i\ne\oo$, $i=0,1$, then we first connect $\C
C_0$ by a simple nodal deformation with any auxiliary non-singular
curve $\C C_0'$ with $\R C_0'\ne\oo$, then deform $\C B_0$ to
obtain $\C B_1$, and finally, connect $\C C_0'$ by a deformation
with $\C C_1$. The case $\R B_0=\R B_1=\oo$ is analogous. \qed

\enddocument